\documentclass{article}
\usepackage{amssymb}
\usepackage{amsmath}
\usepackage[english]{babel}
\usepackage{graphicx}
\usepackage{theorem}

\newtheorem{theorem}{Theorem}

\newtheorem{definition}[theorem]{Definition}
\newtheorem{example}[theorem]{Example}
\newtheorem{lemma}[theorem]{Lemma}

\begin{document}
\title{We can hear (some of) the shape of dented horns}
\author{
Nils Rautenberg \\Ruhr-Universit\"at Bochum, Fakult\"at f\"ur
Mathematik,\\ Universit\"atsstr. 150, D-44780 Bochum, Germany.\\
Email : Nils.Rautenberg@ruhr-uni-bochum.de}
\maketitle

\begin{abstract}
\noindent In this article we construct a family of domains $\Omega \subset \mathbb{R}^2$ with infinite volume such that the Dirichlet Laplacian $\Delta^D$ has purely discrete spectrum and give precise spectral asymptotics for the eigenvalue counting function in terms of the geometry of $\Omega$. This generalizes the well-known asymptotic formula of Hermann Weyl to this class of infinite volume domains. The construction is elementary, uses only the bracketing technique invented by Weyl himself and it is extendable to arbitrary dimensions.
\end{abstract}
\section{Introduction and main result}
One of the seminal contributions to the spectral theory of the Laplacian $\Delta$ is the asymptotic formula given by H. Weyl in 1911 \cite{Weyl} and its many subsequent generalizations. The original result was concerned with the Laplacian $\Delta=-~\mathrm{div }~\mathrm{grad}$ with Dirichlet conditions acting on $L^2(\Omega)$, where $\Omega\subset \mathbb{R}^n$ was some open connected set with piecewise smooth boundary, or domain for short. It made two statements:\\

\noindent First, if the domain $\Omega$ is bounded, then the spectrum is discrete, that is it consists of a sequence of real eigenvalues $0<\lambda_1(\Omega)\leq\hdots \leq \lambda_k(\Omega) \leq \hdots < \infty$, where in this sequence each eigenvalue is counted as many times as its' multiplicity, and the sequence has no accumulation points other than $\infty$.\\

\noindent Second, if the domain is bounded, then not only does this sequence exist, but it shows also a strong universal behaviour guided by the geometry of the domain. If we define
$$N^D_\Omega(E)=\#\{ j : \lambda_j(\Omega)\leq E \}$$ 
to be the counting function of eigenvalues of the Dirichlet Laplacian $\Delta^D$, then the following asymptotic law is true:
$$ N_\Omega^D(E)= \frac{\omega_{n-1}}{(2\pi)^n} \mathrm{vol}_n (\Omega) E^{n/2}+ o(E),$$
where $\mathrm{vol}_n(A)$ denotes the volume of a set measured with the $n$ dimensional Lebesgue measure and $\omega_n$ is a shorthand for $\mathrm{vol}_n(S^n)$, the surface area of the sphere. Weyl also conjectured the existence of a correction term of lower order, encoding the length of the boundary, namely that if the boundary $\partial \Omega$ is piecewise smooth we should get the correction:
$$ N_\Omega^D(E)= \frac{\omega_{n-1}}{(2\pi)^n} \mathrm{vol}_n (\Omega) E^{n/2} - \frac{1}{4} \frac{\omega_{n-2}}{(2\pi)^{n-1}}  \mathrm{vol}_{n-1} (\partial \Omega) E^{(n-1)/2}+o(E^{(n-1)/2}).$$
This conjecture has been proven under the extra assumption that the periodic orbits of the billiard flow on $\Omega$  have measure zero by V. Ivrii in 1980 \cite{Ivr}. Even before the conjecture was proven, M. Kac asked the famous question: "Can one hear the shape of a drum?" \cite{Kac} If the answer was yes it would imply the existence of many more correction terms that ultimately should describe a domain $\Omega$ up to isometry. Sadly, the answer to this question is no, see \cite{Bu}, but it remains open for domains with smooth boundary and it is not clear how large isospectral families of domains can be, or in other words how negative the answer actually is.\\
In addition to this the Weyl law leaves an interesting gap: despite the appealing simplicity of the leading order
asymptotics in these cases, neither boundedness nor finite volume of $\Omega$ are
necessary conditions for purely discrete spectrum of the Dirichlet
Laplacian, denoted from now on by $\Delta^D$. A class of
$2$-dimensional examples for infinite volume domains with discrete
spectrum, as well as their leading order eigenvalue asymptotics
was given by B. Simon in \cite{Sim2}. He considered 2-dimensional
domains of the form:
\begin{equation*}
\Omega_\alpha:=\{(x,y)\in\mathbb{R}^2 : |x|^\alpha|y|< 1\}
\end{equation*}
with $\alpha>0$ and derived the asymptotics:
\begin{equation*}
N^D_\Omega(E)=\begin{cases} \zeta(\alpha)
(\frac{\pi}{2})^{-\alpha}\frac{\Gamma(\frac{1}{2}\alpha+1)}
{\sqrt{\pi}\Gamma(\frac{1}{2}\alpha+\frac{3}{2})}E^{\frac{\alpha+1}{2}}
+ \mathrm{o}(E^\frac{\alpha+1}{2}), \quad \alpha>1\\
\frac{1}{\pi}E\ln{E}+\mathrm{o}(E\log(E)), \quad \alpha=1.
\end{cases}
\end{equation*}
\noindent For $1>\alpha>0$ the first formula holds if one replaces
$\alpha$ by $\alpha^{-1}$. In his paper, Simon raised the question whether this asymptotic behavior and in particular the leading order constant still has a geometric, or physical, interpretation. We will give a partial answer to this question in what follows. We cannot explain the leading terms of Simon's example geometrically, but we will construct a family of infinite volume domains for which such an explanation is possible. However the geometry will not only determine the constant but also the asymptotic growth rate. The methods used are direct and elementary and do not allow to generalize to arbitrary domains immediately, but the results on this examples strongly indicate that additional geometric contributions to the counting function beside the ones seen above might lead to such a law. We begin by defining our domains.

\begin{definition}
\noindent Consider a union of rectangles:
$$\bar{\Omega}=\bar{\Omega}(a,b)=\bigcup\limits_{k=1}^\infty Q_k\subset \mathbb{R}^2$$

\noindent where $Q_k=[a_k,a_{k+1}]\times[0,b_k]$ and the sequences of real numbers $\{a_k\}$ and $\{b_k\}$ satisfy the following properties:
\begin{itemize} 
\item[i)] $b=\{b_k\}$ is monotone, positive and converges to $0$.
\item[ii)] $a=\{a_k\}$ is defined recursively by $a_1=0$, $a_{k+1}=a_k+f(k)$, where $f:\mathbb{N} \to \mathbb{R}^+$ is any function satisfying $f\geq c>0$ for some $c$.
\end{itemize}
Let $\Omega(a,b)=\bar{\Omega}(a,b)^{int}$ denote the interior of this union. In particular it is connected. We will call such a domain $\Omega(a,b)$ \emph{simple}.
\end{definition}

\noindent Because of the bound $f \geq c>0$, all simple domains are unbounded, and depending on the nature of $f$ resp. $a$ and $b$, will have finite or infinite volume. They are horn-shaped in the sense of \cite{Berg} and thus have discrete spectrum. In \cite{Berg}, M. van den Berg gave the precise first order asymptotics of all horn-shaped domains and thus also of all simple domains $\Omega$. While precise, the formula of van den Berg does not give the asymptotics in terms of geometric data that seems compatible with the Weyl law. We will proceed to do this now.
\begin{definition}
Let $\Omega=\Omega(a,b)$ be a simple domain. Define 
$$Q(E)=\left(\bigcup\limits_{k=1}^{n(E)}Q_k\right)^{int},$$ with $$n(E)=\max \{k \in \mathbb{N} ~|~ b_k \geq \frac{\pi}{\sqrt{E}}\}.$$
We call $Q(E)$ a \emph{spectral core} of $\Omega$ for energy $E$.
\end{definition}

\noindent All spectral cores are connected, open and bounded subsets of $\Omega(a,b)$, and they form a nested family that exhausts the domain. Moreover, their geometric properties fully determine the leading term of the eigenvalue counting function:

\begin{theorem}\label{xt}
Let $\Omega=\Omega(a,b) \subset \mathbb{R}^2$ be a simple domain and let $\{Q(E)\}$ be its' family of spectral cores. Assume that the sequence $b$ is summable. Then the following asymptotic law holds true:
$$N_{\Omega}^D(E)=\frac{1}{4\pi} \mathrm{vol}_2(Q(E))E- \frac{1}{4\pi}\mathrm{vol}_1(\partial(Q(E)))E^{\frac{1}{2}}+G(Q(E))+O(E^{\frac{1}{2}}),$$
with $$G(Q(E))=\sum\limits_{k=1}^{n(E)}G_k(E)\mathrm{vol}_2(Q_k),$$
where in each summand $$G_k(E)=\frac{1}{4}[N^D_{Q_k}(E) \mathrm{vol}_2 (Q_k)-\mathrm{vol}_2(B(0, \frac{\sqrt{E}}{\pi}))]$$
is the area error in approximating the volume of a disc $B(0,\frac{\sqrt{E}}{\pi})$ by the sum of the area of all lattice cells of the lattice generated by the dual rectangle of $Q_k$, namely $[0,f(k)^{-1}]\times[0, b_k^{-1}]$ whose corresponding lattice point lie inside the disc. The sum of the three terms is never zero and always at least of order $O(E)$.
\end{theorem}

\noindent In particular, the leading term is completely determined by the geometric data of the domain $\Omega(a,b)$ and more specifically by the geometry of the family of cores $Q(E)$. This true because all three summands are determined by the values of the sequences $a$ and $b$ and for every energy $E>0$, the only data used are the parts of the sequences $a$ and $b$ that determine the shape of $Q(E)$.\\

\noindent The geometric interpretation of the first two summands is straight forward, but the third term $G(Q(E))$ seems rather intricate and elusive with respect to a simple interpretation. It is closely related to the error function of the Gau\ss ian circle problem. Note that the area errors $G_k(E)$ may have a sign and that cancellations are expected to occur. In particular, the sum of these errors $G(Q(E))$ will always oscillate. It is hard to verify to what order the term contributes in a specific example, however it is easy to verify the general bound 

$$|G(Q(E))|\leq c \mathrm{vol}_1(\partial Q(E)) E^{1/2}$$
and we expect that there exist examples where this bound is sharp and of leading order.  We give the terms $\mathrm{vol}_2(Q(E)), \mathrm{vol}_1(\partial Q(E))$ for an example. This shows that both the volume and surface area term may contribute to the leading asymptotics:

\begin{example}
Let $a_{k+1}=a_k+k^3$ and let $b_k=\frac{1}{k^2}$. Then we have:
\begin{itemize}
\item[1)] $n(E)=\lceil \frac{E^{1/4}}{\sqrt{\pi}}\rceil$
\item[2)] $\mathrm{vol}_2(Q(E))=\frac{n(E)(n(E)+1)}{2}=O(E^{1/2})$
\item[3)] $\mathrm{vol}_1(\partial Q(E))=\frac{n(E)^2(n(E)+1)^2}{4}+O(1)=O(E)$
\end{itemize}
\end{example}

\noindent We conclude the introduction with some remarks. It is not necessary to have a summable sequence $b$ in order to obtain correct leading asymptotics, but the asymptotic analysis becomes easier under this condition. In particular, the error bound $O(\sqrt{E})$ in the Theorem depends on this condition. The proof we present here is elementary and relies on the fact that there exists an explicit comparison problem that is asymptotically exact. In particular, the class of simple domains is not the only one for which such a construction is possible, because it is not the only class for which the comparison problem is exact. This also seems more of a technical limitation and not a fundamental obstacle.\\

\noindent The proof we present here also allows to immediately generalize Theorem \ref{xt} to higher dimensions, we give the full proof in two dimensions version for presentation reasons only. If one can decompose a domain $\Omega$ into cuboids, and can do this in such a way that the intersection of the boundary of these cuboids with $\Omega$, the so called \emph{artificial boundary}, is of finite surface area the following generalization of Theorem \ref{xt} is true with error bound $O(E^{(n-1)/2})$:

\begin{theorem}\label{xt2}
Let $\Omega \subset \mathbb{R}^n$ be a domain such that it decomposes into cuboids with summable artificial boundary area. Then there exists again a family of spectral cores $Q(E)\subset\Omega$ such that the following asymptotic law holds true:
\begin{eqnarray*}
N_{\Omega}^D(E)&=&\frac{\omega_{n-1}}{(2\pi)^n} \mathrm{vol}_n(Q(E))E^{n/2}- \frac{1}{4} \frac{\omega_{n-1}}{(2\pi)^{n-1}}\mathrm{vol}_{n-1}(\partial(Q(E)))E^{(n-1)/2}\\
&&+G(Q(E))+O(E^{(n-1)/2}),
\end{eqnarray*}
where $Q(E)$ is the interior of the union of all rectangles sufficiently thick to support eigenvalues of size $E$ or less, and $G(Q(E))$ is the summed up volume error divided by a dual lattice cell for the lattice problem on each cuboid that is part of $Q(E))$.
\end{theorem}

\noindent Based on our analysis on this domain class we conjecture that the notion of spectral core can be generalized to arbitrary domains with discrete spectrum, and that a generalization of the Weyl law in terms of the geometry of such cores is possible. However, the geometry of domains which do not consist of cuboids might contribute to the leading term with even other geometric terms than those visible in this example class. A positive or negative answer will require additional research on low order correction terms of the counting function of compact domains. Figuratively speaking we should be able to hear some of the shape of infinite volume domains, since we can clearly achieve this in the above class.\\

\noindent We even propose a candidate for the generalized notion of core. H. Donnelly obtained in \cite{Don2} necessary and sufficient conditions for a very large class of domains $\Omega\subset (M,g)$ in Riemannian manifolds to have discrete spectrum with respect to $\Delta^D$. Namely, compactness of all sets of the type $\Omega_\epsilon=\{ p \in \Omega ~|~ d_g(p, \partial \Omega) \geq \epsilon \}$ is necessary and sufficient for a domain in this class to have discrete spectrum. When correctly scaled relative to the asymptotic parameter $E$, they seem to provide asymptotically the same geometric information as the cores $Q(E)$. The correct scaling is universal in this example class. The subsets $\Omega_{\pi/\sqrt{E}}$ and $Q(E)$ have asymptotically equivalent growth rates of volume and boundary length. However, a translation of $G(Q(E))$ to $\Omega_{\pi/\sqrt{E}}$ is missing at this point.

\section{Proof of the main result:}

To prove the Theorem \ref{xt} we recall the Dirichlet-Neumann bracketing technique and apply it to our situation.
\begin{lemma}
Let $\Omega(a,n)$ be any simple domain. Then the following bounds are valid by Dirichlet-Neumann bracketing:
$$\sum\limits_{k=1}^\infty N^D_{Q_k}(E)\leq N^D_{\Omega_{a,b}}(E) \leq  \sum\limits_{k=1}^\infty N^{DN}_{Q_k} (E), $$
where $N^D_{Q_k}(E)$ and $N^{DN}_{Q_k} (E)$ are the counting functions of the Laplacian on the single rectangle $Q_k$  with Dirichlet conditions on all vertices respectively Neumann conditions on all vertices that intersect the interior of $\Omega_{a,b}$ and Dirichlet conditions on all other vertices.
\end{lemma}

\noindent The eigenvalue problems in the upper and lower bound are each completely solved. Recall:

\begin{lemma}
Let $Q=[0,a] \times [0,b]$ be a rectangle. Then the spectrum of the Laplacian $\Delta^D$ with Dirichlet condition is given by:
$$\sigma(\Delta^D)=\{ \frac{\pi^2 l^2}{a^2}+ \frac{\pi^2 m^2}{b^2} ~|~ l,m \in \mathbb{N} \},$$
and the spectrum of the mixed Laplacian $\Delta^{DN}$ with Dirichlet conditions on $[0,a] \times \{0\}$ and $[0,a] \times \{b\}$ and Neumann conditions on the remaining boundary components is given by:
$$\sigma(\Delta^{DN})=\{ \frac{\pi^2 l^2}{a^2}+ \frac{\pi^2 m^2}{b^2} ~|~ l \in \mathbb{N} \cup \{0 \}, m \in \mathbb{N} \}.$$
\end{lemma}

\noindent In particular, the series in the upper and lower bound of Lemma 5 are in fact sums for every energy $E>0$, since the thin rectangles $Q_k$ with $k>n(E)$ have no eigenvalue that is small enough to be counted. The main step in the proof of Theorem \ref{xt} is now immediate. If the sequence $b$ of a simple domain $\Omega_{a,b}$ is summable, then the asymptotic bounds formulated in Lemma 5 are sharp to order $\sqrt{E}$, as verified by explicit evaluation of the two comparison problems.

\paragraph{Proof of Theorem \ref{xt}:}
The error between upper and lower bound on each rectangle $Q_k$ is  
$$N^{DN}_{Q_k}(E)- N^D_{Q_k}(E)=\#\{ m \in \mathbb{N} ~|~ \frac{\pi^2 m^2}{b_k^2} \leq E  \}=\lfloor\frac{\sqrt{E}b_k}{\pi}\rfloor.$$
Assume that the sequence $b$ is summable. Then the total error in the approximation satisfies:
$$\sum\limits_{k=1}^\infty N^{DN}_{Q_k}(E)- \sum\limits_{k=1}^\infty N^D_{Q_k}(E)\leq \frac{\sqrt{E}}{\pi} \sum\limits_{k=1}^\infty b_k=O(\sqrt E)$$
Therefore, we can use either counting function to determine the leading terms, and this is now possible because the problem is reduced to an explicit one. We merely need to carefully express the counting function on rectangles in geometric terms, and then sum the contributions from each rectangle. However, this has to be done with care, since errors that do not matter in the asymptotics of a single rectangle become prominent when summed in this way.\\

\noindent The counting function on a single rectangle with Dirichlet conditions is given by:
$$N^D_{Q_k}(E)=\#\{ l,m \in \mathbb{N} ~|~ \frac{\pi^2 l^2}{f(k)^2}+\frac{\pi^2 m^2}{b_k^2} \leq E  \}=\#\{ l,m \in \mathbb{N} ~|~ \frac{ l^2}{f(k)^2}+\frac{m^2}{b_k^2} \leq \frac{E}{\pi^2}  \}.$$
The traditional method to estimate this number is by reinterpreting it as the number of lattice points with positive coordinates of the lattice generated by $Q^*_k=[0,f(k)^{-1}] \times [0,b_k]$, the dual rectangle of $Q_k$, that fit in the disc $B(0,\frac{\sqrt{E}}{\pi})$. This number usually is approximated by dividing the volume of the disc by the volume of a single lattice cell, and is then modified by several correction terms:
\begin{eqnarray*}
N^D_{Q_k}(E)&=& \frac{1}{4} \mathrm{vol}_2(B(0, \frac{\sqrt{E}}{\pi}))\mathrm{vol}_2(Q_k^*)^{-1}-\frac{1}{2}\frac{\sqrt{E}}{\pi} (f(k)+b_k)\\
&&+G_k(E)\mathrm{vol}_2(Q_k^*)^{-1}+R(E)\\
&=& \frac{1}{4\pi} \mathrm{vol}_2(Q_k) E-\frac{1}{4}\frac{\sqrt{E}}{\pi} \mathrm{vol}_1(\partial Q_k) + G_k(E) \mathrm{vol}_2(Q_k) +O(1) 
\end{eqnarray*}
The first two corrections are familiar. We reduce the value by an amount that corresponds to the lattice points on the axis, and we made some error in our volume estimate, accounted for by $G_k(E)$. The final error term $R(E)$ comes from the possible errors at the intersection of the axes and the boundary of the ball and in the origin. But the number of times we miscounted there is bounded by a constant, namely the number of intersections of the axes with the boundary of the disc plus once at the origin. Therefore, the remaining error is not bigger than a constant.\\ 

\noindent  Now, when summing the contributions from each rectangle we find that the contribution from this last error $R(E)$ is negligible, because the sum of it is at most of order $O(n(E))$, and the assumption that the sequence $b$ is summable also implies a bound on $n(E)$. Namely, if $b$ is a summable sequence there must exist a constant $c_1>0$ such that:

$$\sum\limits_{k=1}^n b_k \leq c_1 \sum\limits_{k=1}^n \frac{1}{k}, \quad \forall n \in \mathbb{N}. $$
From this, using the definition of $n(E)$ directly, we find the bound $n(E)\leq c_2 \sqrt{E}$ for some $c_2>0$. Therefore, the asymptotically relevant terms are found from summing the first three terms on each rectangle:

\begin{eqnarray*}
\sum\limits_{k=1}^\infty N^D_{Q_k}(E)=&&\frac{1}{4\pi}\sum \limits_{k=1}^{n(E)}\mathrm{vol}_2(Q_k) E - \frac{1}{4\pi} \sum\limits_{k=1}^{n(E)} \mathrm{vol}_1 (\partial Q_k) E^{1/2}\\ 
&&+ \sum\limits_{k=1}^{n(E)} G_k(E) \mathrm{vol}_2(Q_k) + \mathrm{O}(E^\frac{1}{2}).
\end{eqnarray*}
All summands but the second are already the same as in Theorem \ref{xt}. Now it is easy to see that because $b$ is summable, we have that 
$$\mathrm{vol}_1(\partial Q(E))= \sum\limits_{k=1}^{n(E)} \mathrm{vol}_1 (\partial Q_k) +O(\sqrt{E})$$
and thus, the Theorem is proven. $\hfill\square$\\

\noindent Theorem \ref{xt2} follows by applying the same arguments in higher dimensions. 
\medskip

\noindent \textbf{Acknowledgments.} The author was supported by the project \textit{SFB-TR12},
\textit{Symmetries and Universality in Mesoscopic Systems} founded
by the DFG.

\end{document}